\input amstex
\documentstyle{amsppt}
\magnification=1200
\vcorrection{-0.2in}
\NoRunningHeads
\NoBlackBoxes
\topmatter
\title Homeomorphism Classification of positively curved
manifolds with almost maximal symmetry rank
\endtitle
\author Fuquan Fang \footnote{Supported by CNPq of Brazil, NSFC
Grant 19741002, RFDP
and Qiu-Shi Foundation of China
\hfill{$\,$}}\,\& \,Xiaochun Rong\footnote{Supported partially by
NSF Grant DMS 0203164 and a research grant from Capital normal 
university.
\hfill{$\,$}}
\endauthor
\address Nankai Institute of Mathematics, Nankai University, Tianjing
300071, P.R.C.
\newline
.\hskip3mm Instituto de Matematica, Universidade Federal do Fluminense,
Niteroi, RJ, Brasil
\endaddress
\email ffang\@nankai.edu.cn \hskip4mm fuquan\@mat.uff.br
\endemail
\address Mathematics Department, Rutgers University, New Brunswick,
NJ 08903, U.S.A
\newline
.\hskip3mm Mathematics Department, Capital Normal University,
Beijing, P.R.C.
\endaddress
\email rong\@math.rutgers.edu
\endemail
\abstract We show that a closed simply connected $8$-manifold 
($9$-manifold) of positive sectional curvature on which a 
$3$-torus ($4$-torus) acts isometrically
is homeomorphic to a sphere, a complex projective space or 
a quaternionic projective plane (sphere). We show that
a closed simply connected $2m$-manifold ($m\ge 5$) of 
positive sectional curvature on which 
a $(m-1)$-torus acts isometrically is homeomorphic to a complex 
projective space if and only if its Euler characteristic is 
not $2$. By [Wi], these results imply a homeomorphism classification
for positively curved $n$-manifolds ($n\ge 8$) of almost maximal 
symmetry rank $[\frac{n-1}2]$.  
\endabstract
\endtopmatter
\document

\vskip2mm

\head 0. Introduction
\endhead

\vskip4mm

An interesting problem in Riemannian geometry is the classification
of the positively curved manifolds with large symmetry rank i.e. the
rank of a maximal torus of the isometry group (cf. [Gro]). This
study has gained considerable progress in the recent years
(cf. [FMR], [FR2,3], [GS], [Hi], [HK], [PS], [Ro1], [Wi], [Ya], etc).

Grove-Searle proved that any closed positively curved $n$-manifold 
$M$ has symmetry rank $\le [\frac{n+1}2]$ (the integer part) and 
``$=$'' implies that $M$ is {\it diffeomorphic} to a sphere, 
a lens space or a complex projective space ([GS]).

Very recently, Wilking made a remarkable discovery
that for a closed $n$-manifold of positive sectional curvature
$M$, its a closed totally geodesic $m$-submanifold $N$ (if any)
may capture the homotopy information i.e. $\pi_i(M,N)=0$
for $i\le 2m-n+1$ ([Wi], cf. [FMR]).
Using this result, he proved the following almost $1/2$-maximal rank
theorem: For $n\ge 10$, if a closed simply
connected positively curved $n$-manifold $M$ with symmetry rank
at least $(\frac n4+1)$, then $M$ is {\it homeomorphic} to a sphere or
a quaternionic projective space or {\it homotopically equivalent}
to a complex projective space (for non-simply connected
case, see [Ro1], [Wi])

Subsequently, an analog of the above results for positively
curved manifolds which admit isometric actions by
elementary $p$-groups of large rank has been
obtained in [FR2].

We now state the main results of this paper.

\vskip2mm

\proclaim{Theorem A}

A closed $8$-manifold ($9$-manifold) of positive sectional curvature 
on which a $3$-torus ($4$-torus) acts isometrically is homeomorphic 
to a sphere or a complex projective space or a quaternionic projective plane
$\Bbb HP^2$ (a sphere).
\endproclaim

\vskip2mm

Theorem A implies an extension of the Wilking's almost $1/2$-maximal
rank theorem to dimensions $8$ and $9$. It also arises a natural
question: Is there a {\it critical} symmetry rank between
$\frac n4+1$ and $[\frac{n+1}2]$ so that above which there 
is a homeomorphism classification?

The following provides a partial positive answer.

\vskip2mm

\proclaim{Theorem B}

A closed simply connected $2n$-manifold ($n\ge 5$) $M$ of 
positive sectional curvature on which a $(n-1)$-torus 
acts isometrically is homeomorphic to a complex projective 
space if and only if the Euler characteristic of $M$ is 
not equal to $2$.
\endproclaim

\vskip2mm

As a consequence of Theorems A, B and the Wilking's
almost $1/2$-maximal rank theorem (also, [FR2]), we obtain the following
homeomorphism classification for positively curved manifolds
with almost maximal rank.

\vskip2mm

\proclaim{Corollary C (Almost Maximal Rank)}

For $n\ge 8$, a closed simply connected $n$-manifold $M$ of
positive sectional curvature and almost maximal
symmetry rank $[\frac{n-1}2]$ is homeomorphic to a sphere, a
complex projective space or the quaternionic projective plane
$\Bbb HP^2$.
\endproclaim

\vskip2mm

Corollary C fits nicely between the maximal rank theorem of
Grove-Searle and the almost $1/2$-maximal rank theorem of Wilking
mentioned earlier; in terms of the largeness on symmetry
ranks and different types of classifications: diffeomorphism,
homeomorphism and homotopy. Note that Corollary C is also
sharp in the sense that a sphere, a complex projective space
and a quaternionic projection plane all admit metrics
satisfying Corollary C.

Corollary C is indeed valid for $n=4$ ([HK]), and $n=5$ ([Ro1])
where manifolds in Corollary C were studied. Hence, Corollary C
also generalizes these early results.

To the contrast, the statement of Corollary C no longer holds
in the remaining dimensions $n=6, 7$. There are closed simply
connected $6$- and $7$-manifolds of positive curvature and
almost maximal symmetry rank which are not homeomorphic to any
of rank one symmetric spaces ([AW], [Es]).

Due to the subtleness in dimensions $n=6, 7$, one may like to
take the approach by solving two problems below:

\vskip2mm

\noindent {\bf Problem 0.1}: Show that if $M$ is a positively curved
$6$- (resp. $7$-) manifold of almost maximal rank, then the second
Betti number is at most two (resp. one).

\vskip2mm

\noindent {\bf Problem 0.2}: Classify positively curved $6$- (resp.
$7$-) manifolds of almost maximal rank with $b_2(M)\le 2$ (resp.
$\le 1$).

\vskip2mm

Our last result will provide with a partial solution to Problem 0.2
(for a further development, see [FR3]).

\vskip8mm

\proclaim{Theorem D}

A closed simply connected $6$-manifold ($7$-manifold) of 
positive sectional curvature on which a $2$-torus ($3$-torus) 
acts isometrically is homeomorphic to a sphere if and only if 
its second Betti number vanishes.
\endproclaim

\vskip2mm

Obviously, Theorem D is false if one replaces ``positive curvature''
by ``non-negative curvature''; e.g. the metric product of the
unit spheres, $S^3\times S^3$ and $S^3\times S^4$.
By [GZ], each $S^3$-bundle over $S^4$ admits infinitely many
non-negatively curved metrics. An interesting open question is to
determine if there is a $S^3$-bundle over $S^4$ which is
not diffeomorphic to a standard example and which admits
a metric of positive sectional curvature. Theorem D implies
that such a positively curved metric, if exists, has symmetry
rank at most $2$ (cf. [GKS]).

We now give an indication for the proof of Theorem B.

First, by [Ro1] (see Theorem 1.11) and the almost $1/2$-maximal rank
theorem of Wilking, one concludes that $M$ is homotopy
equivalent to $\Bbb CP^n$. Then, from the
homeomorphism classification of homotopy complex projective
spaces by Sullivan ([Su]), it is easy
to see that $M$ is homeomorphic to $\Bbb CP^n$ if $M$ has
a submanifold $N$ which is homeomorphic to $\Bbb CP^{n-1}$
and the inclusion, $N\hookrightarrow M$, is at least
$3$-connected (see Corollary 1.4).

The proof is to construct a codimension $2$ submanifold
$N$ that meets the above requirements. Observe that combining
the maximal rank theorem and the almost $1/2$-maximal rank
theorem mentioned earlier, it is not hard to see that any invariant
closed
totally geodesic submanifold of codimension $2$ can serve
as the desired $N$ (see Lemma 2.1). Unfortunately, in our
circumstance $M$ may not have any invariant totally geodesic
submanifold of codimension $2$ (see Example 1.8 in [Ro1]).
Hence, the submanifold $N$ we constructed may not be a
totally geodesic submanifold (see Theorem 5.1).

Our strategy is to classify the singular 
structure of an isometric $T^{n-1}$-action. In spirit, this
is similar to [HK] and [Ro1]. The motivation is that the
singular set which is the union of the fixed point sets
of all isotropy groups may provide more topological information
than that from the fixed point set of any single isotropy
group. The two main ingredients in our construction of $N$ are:

\noindent (0.3.1) A classification of the singular set $S$ (as
a stratified space) of the $T^{n-1}$-action (see Theorem 3.6).

\noindent (0.3.2) The orbit space $M^*=M/T^{n-1}$ is homeomorphic
to a sphere (see Theorem 4.1).

The singular set $S$ is stratified, and a stratum is a component
of the singular orbits whose isotropy groups are the same. We will
call the closure of the orbit projection of a singular stratum
a {\it simplex} of $S^*$ (see Section 2), and call the union
of the simplex of dimension $i$ the {\it $i$-skeleton} of $S^*$.

Perhaps it is easier to explain the idea of our construction if
one pretends that the $T^{n-1}$-action is sitting in some isometric
$T^n$-action on $M$ that is invisible to an observer. Then $S^*$
should be the codimension one `skeleton'
of the singular set $(S')^*$ of the $T^n$-action. Observe that
a codimension $2$ invariant totally geodesic submanifold
(which is from the $T^n$-action and which is invisible) in $M$
corresponds
to a top simplex of $(S')^*$. Here (0.3.1) and (0.3.2) enter
to guarantee that one may fill in
a `topological' simplex to replace the top simplex that is
invisible because the codimension one skeleton of the top
simplex is visible in $S^*$ and because whose total space is an 
embedded
sphere in $M^*$ of codimension $3$. Since the embedded sphere is
unknotted ([St]), it must bound a standard topological ball
in $M^*$ (see (0.3.2)). Then, it is not hard to check
that the preimage $N$ of the ball by the orbit projection
has the desired property (Lemma 5.2).

The ingredients in (0.3.1) is the rigidity of the local
singular structure and the connectedness of $S$ (the Frankel's 
theorem). The local rigidity is due to the almost maximality
of symmetry rank and Berger's vanishing theorem (cf. [GS], [Ro1]).
The proof of (0.3.2) is based on (0.3.1) and uses, among
others, the connectedness result of Wilking mentioned earlier.

It is worth to mention that (0.3.1) and (0.3.2) are also
the main ingredients in the proof of Theorems A and D;
using which we may either determine the homology
structure or present $M$ as a union of two standard
pieces from which the homeomorphism type of $M$ is obvious.

For $n=6$ and $7$, $S$ may not be connected and this is
the reason we do not know a classification for $S$ in
the general situation. However, such a classification
is obtained in the case $b_2(M)=0$.

The rest of the paper is organized as follows:

\noindent In Section 1, we collect results that are used in the proofs.

\noindent In Section 2, we prove Theorems A and B for the
case that there is a non-trivial isotropy group with a codimension $2$
fixed point set. Hence, we may assume that in
the rest sections, any isotropy group has fixed point set
codimension at least $4$.

\noindent In Section 3, we classify the singular structure.

\noindent In Section 4, we prove that the orbit space $M^*$ is
homeomorphic to a sphere.

\noindent In Section 5, we prove Theorem B.

\noindent In Section 6, we prove Theorem A.

\noindent In Sections 7 and 8, we prove Theorem D.

\vskip4mm

\head 1. Preliminaries
\endhead

\vskip4mm

In this section, we will collect results that will be served
as tools in the rest of the paper.

\vskip4mm

\subhead a. Homeomorphic classification of a homotopy type
\endsubhead

\vskip4mm

The generalized Poincar\'e conjecture says that any homotopy
$n$-sphere is homeomorphic to a sphere $S^n$. By the famous
works of Smale and Freedman, only the case $n=3$ (the original
Poincar\'e conjecture) remains unsolved.

\vskip2mm

\proclaim{Theorem 1.1 ([Fr])}

Any homotopy $4$-sphere is homeomorphic to $S^4$.
\endproclaim

\vskip2mm

\proclaim{Theorem 1.2 ([Sm])}

For $n\ge 5$, any homotopy $n$-sphere is homeomorphic to $S^n$.
\endproclaim

\vskip2mm

Contrast to being a sphere, the homotopy type of a closed
manifold may contain many homeomorphism types. In [Su],
Sullivan classified the homeomorphism types of a homotopy
complex projective space.

Let $h: M\to \Bbb CP^n$ be a homotopy equivalence. For any $\Bbb CP^i
\subset \Bbb CP^n$, let $h^{-1}(\Bbb CP^i)\subset M$ denote the
transverse submanifold
(the preimage of a map homotopic to $h$ and transverse to
$\Bbb CP^i$.)
Let $\sigma _h(i)$ be the signature (an integer) of $h^{-1}(\Bbb
CP^i)$
if $i$ is even (an integer), and the Kervaire invariant (a mod $2$
integer) if $i$ is odd (cf. [Su]). The invariant $\sigma _h(i)$
depends only on the homotopy class $h$. Indeed, if $N^{2i}\subset M$
is another manifold homologous to $h^{-1}(\Bbb CP^i)$, then
$\sigma _h(i)$ is
equal to the signature  (resp. Kervaire invariant) of $N^{2i}$.

\vskip2mm

\proclaim{Theorem 1.3 ([Su])}

Let $hS(\Bbb CP^n)$ denote the set of homeomorphism classes of
closed manifolds homotopy equivalent to $\Bbb CP^n$. For $i=1,
\cdots, [n/2]$,
$\sigma _h(i)$ defines an one-to-one correspondence between
$hS(\Bbb CP^n)$ and the set $\prod_{i=1}^{[n/4]} (\Bbb Z
\times\Bbb Z_2)$.
\endproclaim

\vskip2mm

Recall that the Kervaire invariant is zero if a manifold of
dimension $(4i+2)$ has no middle dimensional homology, e.g.
manifolds homotopy equivalent to $\Bbb C P^{2i+1}$. On the
other hand, the signature is a homotopy invariant.  As an
immediate corollary of Theorem 1.3, we conclude

\vskip2mm

\proclaim{Corollary 1.4}

Assume that $M\in hS(\Bbb CP^n)$ has a codimension $2$
submanifold $N$ which is homeomorphic to $\Bbb CP^{n-1}$.
Then $M$ is homeomorphic to $\Bbb CP^n$ under one of
the following conditions:

\noindent (1.4.1) $N$ represents a generator for $H_{2n-2}(M)$.

\noindent (1.4.2) The inclusion, $N\hookrightarrow M$, is at
least $3$-connected.
\endproclaim

\vskip4mm

\subhead b. $G$-spaces
\endsubhead

\vskip4mm

Let $G$ denote a compact Lie group. The following two results,
which concern a $G$-space, will be frequently quoted in our
proofs. The first one is the so-called homotopy
lifting property (cf. [Br]).

\vskip2mm

\proclaim{Lemma 1.5 (Homotopy Lifting)}

Let $G$ be a compact Lie group, and let $M$ be a connected
$G$-manifold. If there is a connected $G$-orbit,
then $p_*: \pi_1(M)\to \pi_1(M/G)$ is surjective.
\endproclaim

\vskip2mm

Let $F(G,M)$ denote the $G$-fixed point set. It is well known
that if $G$ is abelian, then topology of $F(G,M)$ is closely
related to the topology of $M$ (cf. [Hs]). Note that given
an $G$-invariant metric, a component of $F(G,M)$ is a totally
geodesic submanifold of even codimension.

\vskip2mm

\proclaim{Theorem 1.6 (Fixed point)}

Let $M$ be a closed $G$-space.

\noindent (1.6.1) If $G=T^k$ a torus, then the Euler characteristic
$\chi(M)=\chi(F(T^k,M))$.

\noindent (1.6.2) If $G=\Bbb Z_p^k$ ($p$ is a prime), then $\chi(M)=
\chi(F(\Bbb Z_p^k,M))\,(\text{mod}\,\, p)$
\endproclaim

\vskip2mm

\vskip4mm

\subhead c. Positive curvature with large symmetry rank
\endsubhead

\vskip4mm

The rank of a Lie group is the rank of a maximal torus.
Hence, $\text{Symrank}(M)=k$ implies that $M$ admits an
isometry torus $T^k$-action. When a $T^k$-invariant
metric has positive sectional curvature, the basic fact
is the following Berger's vanishing theorem ([Ko], [Ro2]):

\vskip2mm

\proclaim{Theorem 1.7}

Let $M$ be a closed $n$-manifold
of positive sectional curvature. Assume that $M$ admits
an isometric $T^k$-action.

\noindent (1.7.1) If $n$ is even, then the fixed point set
is not empty.

\noindent (1.7.2) If $n$ is odd, then there is a circle orbit.
\endproclaim

\vskip2mm

A consequence of Theorem 1.7 is that a large symmetry rank
implies a totally geodesic submanifold of small codimension.
Recently, Wilking [Wi] discovered that a totally geodesic submanifold
of small codimension may capture a lot homotopy information
of $M$ (see [FMR] for a further development).

A map $f: N\to M$ is called $(i+1)$-connected, if $f_*:
\pi_q(N)\to \pi_q(M)$ is an isomorphism for $q\le i$ and
onto for $q=i+1$.

\vskip2mm

\proclaim{Theorem 1.8 ([Wi])}

Let $M$ be a closed $n$-manifold of positive
sectional curvature, and let $N$ be a closed totally geodesic
$k$-submanifold. If there is a Lie group $G$ that acts
isometrically on $M$ and fixes $N$ pointwisely, then the
inclusion map is $(2k-n+1+C(G))$-connected, where $C(G)$ is
the dimension of a principal orbit of $G$.
\endproclaim

\vskip2mm

In the extremal case when a $T^k$-fixed point set has the maximal
possible dimension (e.g. a circle group has a fixed point set
codimension $2$ or a torus has a fixed point set codimension
$4$), one has the Grove-Searle diffeomorphism classification
for the maximal symmetry rank.

\vskip2mm

\proclaim{Theorem 1.9 ([GS])}

Let $M$ be a closed $n$-manifold of positive sectional curvature.

\noindent (1.9.1) If $M$ admits an isometric circle action
with a fixed point set of codimension $2$, then $M$ is diffeomorphic
to a lens space $S^n/\Bbb Z_p$ ($p$ may be $1$) or a complex
projective space.

\noindent (1.9.2) If $M$ admits an isometric $T^k$-action,
then $k\le [\frac{n+1}2]$ and ``$=$'' implies a circle
subgroup whose fixed point set has codimension $2$.
\endproclaim

\vskip2mm

The above Theorems 1.8 and 1.9 play an important role in
the proof of the following almost $1/2$-maximal rank theorem
of Wilking.

\vskip2mm

\proclaim{Theorem 1.10 ([Wi])}

For $n\ge 10$, let $M$ be a closed simply connected $n$-manifold
of positive sectional curvature. If $\text{Symrank}(M)\ge \frac n4+1$,
then $M$ is homeomorphic to a sphere or a quaternionic projective
space or $M$ is homotopically equivalent to a complex projective
space.
\endproclaim

\vskip2mm

The following theorem in [Ro1] on the singular structure
of a positively curved manifold with almost maximal
symmetry rank will be used.

\vskip2mm

\proclaim{Theorem 1.11 ([Ro1])}

For $n\ge 8$, let $M$ be a closed simply connected $n$-manifold
of positive sectional curvature. Assume that $M$ admits an
isometric $T^{[\frac{n-1}2]}$-action such that there is no
circle subgroup
with fixed point set of codimension $2$.

\noindent (1.11.1) If $n=8$, then the fixed points are isolated
whose number is $2, 3$ or $5$.

\noindent (1.11.2) If $n=2m>8$, then the fixed points are isolated
whose number is either $2$ or $(n+1)$.

\noindent (1.11.3) If $n=2m+1$, then the circle orbits are isolated
whose number is $(m+1)$.
\endproclaim

\vskip2mm

A consequence of Theorems 1.11 and 1.6 is

\vskip2mm

\proclaim{Corollary 1.12}

For $n\ge 4$, let $M$ be a closed $2n$-manifold of positive
sectional curvature and almost maximal symmetry rank $(n-1)$.
Then the Euler characteristic of $M$, $\chi (M)=2$ or $(n+1)$
or $3$ when $n=4$. In particular, $M$ is not homotopically equivalent
to $\Bbb HP^n$ for $n\ne 2$.
\endproclaim

\vskip4mm

\head 2. A Special Case of Theorems A and B
\endhead

\vskip4mm

In this section, we prove the case of Theorems A and B that there
is a non-trivial isotropy group with codimension $2$ fixed point set
(Lemma 2.1). We then show that in the complementary situation,
all isotropy groups are connected (Lemma 2.3). This condition is
crucial to obtain a complete classification of the singular set
(see Section 3).

\vskip2mm

\proclaim{Lemma 2.1}

Let $M$ be a closed simply connected $n$-manifold of positive
sectional curvature. Assume that $M$ admits an isometric
$T^{[\frac{n-1}2]}$-action. If there
is a non-trivial isotropy group with codimension $2$ fixed point set,
then $M$ is homeomorphic to a sphere or a complex projective space.
\endproclaim

\vskip2mm

\demo{Proof} First, we may assume $n\ge 6$ since the cases of $n=4$ and
$n=5$ are covered by [HK] and [Ro1].

Let $H$ denote an isotropy group with a codimension $2$ fixed point
component $F_0\subset M$. Since $H$ acts effectively on a normal
$2$-disk of $F_0$, $H=S^1$ or $H$ is a finite cyclic group. By
Theorem 1.9, we may assume that $H$ is finite. Since
$T^{[\frac{n-1}2]}/H$ acts effectively on $F_0$, $F_0$
has the maximal symmetric rank. Since $M$ is simply connected,
$F_0$ is simply connected (Theorem 1.8), and since the induced
metric has positive sectional curvature, $F_0$ is diffeomorphic
to $S^{n-2}$ or $\Bbb CP^{\frac{n-2}2}$ (Theorem 1.9).

Case 1. If $F_0\overset{\text{diff}}\to \simeq S^{n-2}$, since
$F_0\hookrightarrow M$
is $(n-3)$-connected, by the Hurewicz theorem and the
Poincar\'e duality it is easy to see that $M$ is a homotopy
sphere. By Theorem 1.2, $M$ is homeomorphic to $S^n$.

Case 2. If $F_0\simeq \Bbb CP^{\frac{n-2}2}$. Since $F_0
\hookrightarrow M$ is at least $3$-connected, by Corollary 1.4
$M$ is homeomorphic to $\Bbb CP^\frac n2$.
\qed\enddemo

\vskip2mm

The following are cases where Lemma 2.1 may apply.

\vskip2mm

\proclaim{Lemma 2.2}

Let $M$ be a closed $n$-manifold of positive sectional curvature.
Assume that $M$ admits an isometric $T^{[\frac{n-1}2]}$-action.
Then there is an isotropy group with fixed point set of
codimension $2$ if one of the following conditions holds:

\noindent (2.2.1) For $n=2m$, the fixed point set
$F(T^{[\frac{n-1}2]},M)$ has dimension $>0$.

\noindent (2.2.2) For $n=2m+1$, there are non-isolated circle orbits
(e.g. if the fixed point set is not empty).

\noindent (2.2.3) There is a non-trivial finite isotropy group.
\endproclaim

\vskip2mm

\demo{Proof} (2.2.1) Let $F_0$ denote a component of
$F(T^{[\frac{n-1}2]},M)$ of positive dimension.
Since $T^{[\frac{n-1}2]}$ acts effectively on the normal
space $T^\perp_xF_0$ via differentials, $T^\perp_xF_0\simeq
\Bbb R^{n-2}$ and thus $T^{[\frac{n-1}2]}$ has a circle
subgroup $S^1$ with fixed point set codimension $2$ in
$T^\perp_xF_0$. Consequently, $\dim(F(S^1,M))=n-2$.

(2.2.2) Let $H\subset T^{[\frac{n-1}2]}$ be an isotropy group
of a non-isolated circle orbit, and let $F_0$ be a component of
the fixed point set $F(H,M)$ of dimension $2i+1>1$. Then
$H$ has rank $[\frac {n-3}2]$. Since $H$ acts effectively
on the normal space of $F_0$, as in (2.2.1) we conclude
the desired result.

(2.2.3) Let $x\in M$ such that the isotropy group $H$ at $x$ is
finite. Let $F_0$ denote a $H$-fixed point component containing $x$.
Then $T^{[\frac{n-1}2]}\cong T^{[\frac{n-1}2]}/H$ acts effectively
on $F_0$. Since $F_0$ is a closed totally geodesic submanifold of
even codimension, by Theorem 1.9 it follows that
$\dim(F_0)=n-2$.
\qed\enddemo

\vskip2mm

If (2.2.3) does not occur, then we have

\vskip2mm

\proclaim{Lemma 2.3}

Let $M$ be a closed $n$-manifold with a $T^k$-action. Then there is no
finite isotropy group if and only if every isotropy group is connected.
\endproclaim
\vskip2mm

\demo{Proof} It suffices to prove the necessity. We argue by contradiction.
Assume that
$T^\ell\times A$ is an isotropy group of some $x\in M$, where $\ell >0$
is minimal and $A\ne 1$ is a finite abelian group. By the slice theorem
([Br]), the orbit $T^k(x)$ has a
tubular neighborhood $U=T^k\times_{T^\ell\times A}D^\perp_\epsilon$, and
the isotropy groups of the $T^k$-action in $U$ are in one-to-one
correspondence to the isotropy subgroups of $T^\ell\times A$-action in
$D^\perp_\epsilon$, where $D^\perp_\epsilon$ is the
$\epsilon$-ball in the normal slice of $T^k(x)$. Hence, it suffices to
prove that the linear $T^\ell\times A$-action on $\partial D^\perp_\epsilon$
has a non-trivial finite isotropy group.

Assume that $\Bbb Z_p$ is a cyclic subgroup of $A$ of prime order. Consider
the action of $\Bbb Z_p^\ell\times \Bbb Z_p$ on $\partial D^\perp_\epsilon$,
where $\Bbb Z_p^\ell$ is the subgroup of $T^\ell$.  By [Br] III. Theorem
10.12 (and its remark) one sees that, there are at least two isotropy
groups of $p$-rank $\ell$ in $\Bbb Z_p^\ell\times \Bbb Z_p$. Therefore,
there is at least one $p$-rank $\ell$ subgroup of $T^\ell \times A$
but not contained in $T^\ell$, which acts on $\partial D^\perp_\epsilon$
with non-empty fixed point set. Note that every isotropy group of this
fixed point set has the form $T^s\times B$, where $s\le \ell -1$ and $B
\cong \Bbb Z_p^{\ell-s}$. A contradiction to the minimality of $\ell$.
\qed\enddemo

\vskip4mm

\head 3. Singular Structures in The Reduced Case
\endhead

\vskip4mm

By Lemma 2.1, the proofs of Theorems A, B and D may be reduced to the
case where there is no non-trivial isotropy group with codimension
$2$ fixed point set. By Lemmas 2.2 and 2.3, we may assume the action
of $T^{[\frac {n-1}2]}$ satisfies that

\noindent (3.1.1) All isotropy groups are connected.

\noindent (3.1.2) The union of singular orbits is of codimension $4$
(cf. [Ko]).

Note that by Lemma 2.2, (3.1.2) implies the following:

\noindent (3.1.2a) If $n=2m$, then all fixed points are isolated.

\noindent (3.1.2b) If $n=2m+1$, then all circle orbits are isolated.

For an example of a torus action satisfying the above, see
Example 1.8 in [Ro1].

The goal of this section is to classify the singular set of an
isometric $T^{[\frac {n-1}2]}$-action on a positively curved
$n$-manifold $M$ which satisfies (3.1.1) and (3.1.2) (see Theorem  3.6).
This result will be used in the proofs of Theorem 4.1 and Theorems
A, B and D (see the discussion following (0.3.1)).

Let $S$ denote the singular set i.e. the union of non-principle
orbits. For an isotropy
group $H$, let $S_{(H)}$ denote the union of orbits whose
isotropy group is $H$. Then $S$ is stratified by singular
strata i.e. components of $S_{(H)}$ for all isotropy group
$H\ne 1$. The $i$-th skeleton of $S$, $S^{(i)}$, is the union of
singular strata whose dimensions are at most $i$. Let
$p: M\to M^*$ denote the orbit projection. We call $S^*=p(S)$
the singular set of $M^*$. Clearly, the stratification of
$S$ descends to a stratification on $S^*$.

Note that the closure of a singular stratum is a fixed point
component of $H$. We will call a component of the closure
$\bar S_{(H)}$ a {\it simplex}.

Let $\Delta_s$ denote the standard $s$-simplex in $\Bbb R^{s+1}$
whose vertices are $(1, 0, \cdots ,0, 0)$,
$(0, 1, 0, \cdots, 0)$, $\cdots$, $(0, 0, \cdots, 0, 1)$.
Then $\Delta_s$ is a  stratified space with the natural
stratification: a vertex is a $0$-stratum, an open
$1$-simplex (edge) is a $1$-stratum, ...., an open
$s$-simplex is an $s$-stratum. Let $\Delta_s^k$ denote
the $k$-skeleton of $\Delta_s$.

Two stratified spaces $S^*_1$ and $S^*_2$ are isomorphic,
if there is a homeomorphism $f: S^*_1\to S^*_2$ which maps
a stratum onto a stratum.

\vskip4mm

\subhead a. The singular structure of the maximal rank
\endsubhead

\vskip2mm

We will see that the singular set $S^*$ of an almost
maximal symmetry rank action is the codimension one skeleton of a
singular set in some standard example of the maximal
symmetry rank; for which we will now describe its
singular structure.

It is an easy exercise to verify the following three model
cases (Lemmas 3.2-3.4).

\vskip2mm

\proclaim{Lemma 3.2}

Let $T^{m+1}$ denote a maximal torus of $O(2m+2)$. Then $T^{m+1}$
acts isometrically on the unit sphere $S^{2m+1}$ with the
orbit space $S^{2m+1}/T^{m+1}$ homeomorphic to $\Delta_m$ and the
singular set $S^*$ is
isomorphic to $\partial\Delta_m=\Delta_m^{m-1}$ as stratified
space. Moreover, a simplex of $S$ is a totally geodesic sphere
and its orbit projection is a simplex of $\Delta_m$.
\endproclaim

\vskip2mm

\proclaim{Lemma 3.3}

Let $(\Bbb CP^m,T^m)$ denote the complex projective space with
the Fubini-Study metric and
$T^m$ be a maximal torus of $\text{Isom}(\Bbb CP^m)$. Then
the orbit space $\Bbb CP^m/T^m$
is homeomorphic to $\Delta_m$ and the singular set $S^*$ is
isomorphic to $\partial\Delta_m=\Delta_m^{m-1}$ as stratified
space. Moreover, a simplex of $S$ is a totally geodesic
complex projective space and its orbit projection is a
simplex of $\Delta_m$.
\endproclaim

\vskip2mm

We now consider the standard $(S^{2m},T^m)$. Observe that since
there are only
two isolated $T^m$-fixed points, $S^*$ with its simplices defined
in the above does not form a simplicial complex. Because of this,
one may like to refine the stratification of $S^*$ to a simplicial
complex structure as follows.

Let $S^{2m-1}$ denote the subsphere in $S^{2m}$ consisting
of points whose distance from a $T^{2m}$-fixed point is
$\pi/2$. Then the $T^m$-action preserves $S^{2m-1}$. Let
$S^*_1$ denote the singular set of $S^{2m-1}/T^m$.
Note that each $i$-simplex $B_i$ of $S^*_1$ is contained
in a unique $(i+1)$-simplex $A_{i+1}$ in $S^*$. We
then divide $A_{i+1}$ into two new $(i+1)$-simplices
whose intersection is $B_i$. These new simplices also
form a stratification of $S^*$, denoted by
$\text{rf}(S^*)$.

Let $\Sigma(\Delta_{m-1}^{m-2})$ denote the two sides
suspension of $\Delta_{m-1}^{m-2}$ obtained by
identify the two ends of $\Delta_{m-1}^{m-2}\times [-1,1]$
with two points. Then $\Sigma(\Delta_{m-1}^{m-2})$
has a natural stratification which also forms a simplicial
complex.

\vskip2mm

\proclaim{Lemma 3.4}

Let $T^m\subset O(2m+1)$ denote a maximal torus.
Then $S^{2m}/T^m$ is homeomorphic to $\Delta_m$ and $\text{rf}(S^*)$
is isomorphic to $\Sigma(\Delta_{m-1}^{m-2})$ as
stratified space. Moreover, a simplex of $S$ is a totally geodesic
sphere and the orbit projection preserves vertices and maps
a simplex of $S$ to the suspension of a simplex of
$\Sigma(\Delta_{m-1}^{m-2})$.
\endproclaim

\vskip2mm

\proclaim{Lemma 3.5}

Let $M$ be a simply connected positively curved closed $n$-manifold
with an isometric $T^{[\frac{n+1}2]}$-action. Then the orbit space
$M^*$ is homeomorphic to $\Delta_m$, where $n=2m$ or $2m+1$.
Moreover,

\noindent (3.5.1) If $n=2m+1$, then $S\cong S(S^{2m+1},T^{m+1})$
and $S^*\cong \partial \Delta_m$.

\noindent (3.5.2) If $n=2m$ and $\chi(M)=2$, then
$S\cong S(S^{2m},T^m)$ and $S^*\cong
\Sigma (\Delta_{m-1} ^{m-2})$.

\noindent (3.5.3) If $n=2m$ and $\chi(M)=m+1$, then
$S\cong S(\Bbb CP^m,T^m)$ and
$S^*\cong \partial \Delta_m$.
\endproclaim

\vskip2mm

\demo{Proof} Since the proofs are similar, we will only present
a proof for (3.5.1). First, there are $(m+1)$ isolated circle
orbits or equivalently, $S^*$ has $(m+1)$-vertices
(if $n\ge 8$, see Theorem 1.11; and the case
for $n\le 7$ can be easily verified). By Lemma 3.8
below, around each fixed vertex $v$ there are $m$
$1$-simplices whose end ($\ne v$) are all distinct.
The maximality implies that every $(i+1)$-vertices
determines a $i$-simplex of $S^*$ and this clearly
implies the desired result.
\qed\enddemo

\vskip4mm

\subhead b. The singular structure of almost maximal rank
\endsubhead

\vskip4mm

The goal of this subsection is to prove the following
classification result:

\vskip2mm

\proclaim{Theorem 3.6 (Classification of singularities)}

Let $M$ be a closed simply connected $n$-manifold of positive
sectional curvature. Assume that $M$ admits an isometric
$T^{[\frac{n-1}2]}$-action satisfying (3.1.1) and (3.1.2).

\noindent (3.6.1) If $n=2m\ge 8$ and $\chi(M)=m+1$, then
$S^*\cong \Delta_{m}^{m-2}$ and $S\cong S^{(2m-4)}(\Bbb CP^m,T^m)$.

\noindent (3.6.2) If $n=2m+1\ge 9$, then $S^*\cong \Delta_{m}^{m-2}$ and
$S\cong S^{(2m-3)} (S^{2m+1},T^m)$.

\noindent (3.6.3) If $n=2m\ge 8$ and $\chi(M)=2$, then
$\text{rf}(S^*)\cong \Sigma (\Delta_{m-1}^{m-3})$ and $S\cong S^{(2m-4)}
(S^{2m},T^m)$.

\noindent In particular, a simplex of $S$ is a totally geodesic
sphere in (3.6.1) and (3.6.3) and a complex projective space
in (3.6.2). In any case, the orbit space of a simplex
is contractible.
\endproclaim

\vskip2mm

Roughly speaking, Theorem 3.6 says that the singular set $S$
(resp. $S^*$) is a codimension $2$ (resp. codimension one) skeleton
of some corresponding model space of maximal rank (see
Lemmas 3.2-3.4).

\vskip2mm

\remark{Remark \rm 3.7} Consider a positively curved $6$- or
$7$-manifold $M$ with (almost maximal) symmetry rank $2$ and
$3$ respectively. Note that unlike the case of $n\ge 8$, the
singular set $S$ of $M$ may not be connected. This is a problem
to classify the global singular structure of $M$ (cf. [FR3]).
\endremark

\vskip2mm

Similar to the situation of the maximal rank, the local
singular structure in our circumstance is rigid.

\vskip2mm

\proclaim{Lemma 3.8 (Local singular structure)}

Let the assumptions be as in Theorem  3.6.

\noindent (3.8.1) If $n=2m$, then around each fixed point
there are $\frac {m!}{i!(m-i)!}$-many $2i$-dimensional singular
strata whose isotropy groups have rank $(m-1-i)$, $1\le i\le m-2$.

\noindent (3.8.2) If $n=2m+1$, then around each isolated circle orbit
there are $\frac {m!}{i!(m-i)!}$-many $(2i+1)$-dimensional
strata whose isotropy groups have rank $(m-1-i)$, $1\le i\le m-2$.

\noindent (3.8.3) If $H$ is an isotropy group of rank $r$, then
$\dim (S_{(H)})=n-2(r+1)$.

\noindent (3.8.4) If $n\ge 8$ then the singular set $S$ is connected.
\endproclaim

\vskip2mm

\demo{Proof} (3.8.1) Consider the isotropy representation of
$T^{m-1}$ at an isolated fixed point. This is clearly a
sub-representation
of the standard linear $T^{m}$-action on $D^{2m}$ (cf. (3.2.1)). Since
(3.1.2), $T^{m-1}$ must transversally intersect every torus subgroup
$\{ (u_1, \cdots , u_m )\in T^m\subset \Bbb C^m| \text{ where only
one } u_i=1 \text{ for } 1\le i\le m\}$. Therefore, by the standard
linear algebra,
around the fixed point there are exactly $\left(\matrix m \\i
\endmatrix\right)$-many $2i$-dimensional strata whose isotropy groups
have rank $m-1-i$, $1\le i\le m-2$.

(3.8.2) Let $\Cal O_x$ be an isolated circle orbit with isotropy group
$H\cong T^ {m-1}$. Consider the isotropy representation of $H$ at the
normal slice $D^ {2m}$ of $\Cal O_x$. Since the isotropy groups of
the linear
$H$-action on the slice are in one-to-one correspondence with the
isotropy
groups of the $T^m$-action on a slice neighborhood of $\Cal O _x$.
The desired result follows from (3.8.1).

(3.8.3) If $n=2m$ (resp. $n=2m+1)$, then $S_{(H)}$ contains a
fixed point (resp. an isolated circle orbit). The desired result
follows
from (3.8.1) and (3.8.2).

(3.8.4) Since $m\ge 4$, by (3.8.1) over each vertex of $S^*$, $S$ has a
codimension $4$ simplex and any two codimension $4$ simplices have a
non-empty intersection (Frankel theorem). This shows that $S$ is connected.

\qed\enddemo

\vskip2mm

\demo{Proof of Theorem  3.6}

Since each simplex of $S$
is a closed totally geodesic submanifold, each
codimension $4$ simplex is simply connected (Theorem 1.8).
Since a codimension $4$ simplex has the maximal symmetry
rank, it is diffeomorphic to a sphere or a complex projective
space.

(3.6.1) Note that $\Delta_m^{m-2}$ is isomorphic to
a simplicial complex with $(m+1)$-vertices such that
every $(m-1)$-vertices spans a simplex isomorphic to $\Delta_{m-2}$.

By Theorem 1.11, $S^*$ has $(m+1)$ vertices. By Lemma 3.8,
at each (fixed) vertex $v$ there are $(m+1)$ $1$-simplices
containing $v$ and any $(m-1)$ of these $1$-simplices
are contained in a (top) $(m-2)$ simplex. By the above discussion,
the corresponding top simplex of $S$ is a sphere
($n$ is odd) with the maximal symmetry rank. This
implies that a top simplex of $S^*$ is isomorphic to
$\Delta_{m-2}$. Since at $v$, each $1$-simplex
contains another vertex (the other end point), we conclude
that every $(m-1)$ vertices spans a simplex $\Delta_{m-2}$.
We then conclude the desired result with the observation
that $S^*$ completely determines $S$ when
specifying the preimage of every top simplex of $S^*$.

(3.6.2) We omit the proof since it is exactly the same
argument if replacing the assertion that a codimension
$4$ simplex of $S$ is diffeomorphic to a sphere by
that a top simplex is diffeomorphic to a complex
projective space.

(3.6.3). In this case, $S$ (and $S^*$) has exactly
two vertices (Theorem 1.11). The vertices correspond to
the two vertices of the suspension $\Sigma (\Delta_{m-1}^{m-3})$.
\qed\enddemo

\vskip2mm

\proclaim{Lemma 3.9}

Let the assumptions be as in Theorem 3.6. If $n=8$ and
$\chi(M)=3$, then $S\cong S(\Bbb HP^2,T^3)$ and $S^*
\cong S^*(\Bbb HP^2,T^3)$, where $\Bbb HP^2$ is the
quaternionic projective space with the standard metric
and $T^3$ is a maximal torus in $\text{Isom}(\Bbb HP^2)$.
\endproclaim

\vskip2mm

\demo{Proof} It is an easy exercise to see that $S^*(\Bbb HP^2,T^3)$
is isomorphic to a $2$-dimensional singular polyhedron with
three vertices $v_1, v_2, v_3$, six segments $\alpha_1,
\alpha _2$ (joining the vertices $v_1, v_2$), $\beta _1,
\beta _2$ (joining the vertices $v_2$ and $v_3$),
$\gamma _1, \gamma _2$ (joining the vertices $v_3$ and $v_1$)
and seven faces with boundaries $\alpha _1\alpha _2$,
$\beta _1\beta _2$, $\gamma _1\gamma _2$, $\alpha _1
\beta_1\gamma _1$, $\alpha _1\beta _2\gamma _2$,
$\alpha _2\beta _2\gamma _1$, $\alpha _2\beta _1\gamma _2$.
The preimage of a face by the orbit projection, which is a
codimension $4$ simplex, is diffeomorphic to a sphere or a
complex projective space depending the number of the
vertices in the simplex.

We now consider $(M^8,T^3)$. By Theorem 1.11, $S^*$ has exact
three vertices. Together
with the fact that $S^*$ is connected and that at each
vertex there are $4$ $1$-simplex imply that the
$1$-skeleton of $S^*$ is isomorphic to the $1$-skeleton
of $S^*(\Bbb HP^2,T^3)$. Since at a fixed vertex $v$ of
$S$ there are six $4$-simplices, this implies that
$S^*$ is isomorphic to $S^*(\Bbb HP^2,T^3)$. Since
each $4$-simplex is diffeomorphic to either $S^4$
or $\Bbb CP^2$, we then conclude the desired result.
\qed\enddemo

\vskip4mm

\head 4. Orbit Spaces of The Reduced Case
\endhead

\vskip4mm

The goal of this section is to prove the following theorem
(see the discussion following (0.3.2)).

\vskip2mm

\proclaim{Theorem 4.1}

Let $M$ be a closed simply connected $n$-manifold of positive
sectional curvature. Assume that $M$ admits an isometric
$T^{[\frac{n-1}2]}$-action satisfying (3.1.1) and (3.1.2). Then
$M^*$ is a closed topological manifold and for $n\ge 8$,
$M^*$ is homeomorphic to a sphere.
\endproclaim

\vskip2mm

We first verify a special case of Theorem 4.1.

\vskip2mm

\proclaim{Lemma 4.2}

Let the assumptions be as in Theorem 4.1. Assume that $M=S^n$ is
the standard unit sphere. Then $M^*$ is homeomorphic to a sphere
for all $n$.
\endproclaim

\demo{Proof} We will only present a proof for the case $n=2m+1$
since the proof for $n=2m$ is similar.

Note that $T^m\subset \text{Isom}(S^{2m+1})=O(2m+2)$
and thus we may assume that $T^m$ is contained in some
maximal torus $T^{m+1}\subset O(2m+2)$. By Lemma 3.2,
$S^{2m+1}/T^{m+1}$ is homeomorphic to $\Delta_m$.
Let $X=S^{2m+1}/T^m$. Then $S^1=T^{m+1}/T^m$ acts on
$X$ such that $X/S^1\simeq \Delta_m$.

Note that a singular $T^m$-orbit is also a singular
$T^{m+1}$-orbit and any singular $T^{m+1}$-orbit
projects to a point in $X$ which is a $S^1$-fixed
point. Hence, the $S^1$-action on $X$ is semi-free
with fixed point set the boundary of $\Delta_{m}$.
i.e. every orbit over the interior of $\Delta_{m}$
is principal. Therefore, $X$ is homeomorphic to
$(\partial \Delta_{m}\times D^2) \cup _\partial
(\Delta_{m}\times S^1)\approx S^{m+1}$.
\qed\enddemo

\vskip2mm

\proclaim{Lemma 4.3}

Let the assumptions be as in Theorem 4.1. Then $M^*$ is a closed
topological manifold.
\endproclaim

\demo{Proof}
It suffices to prove that any {\it maximal singular point}
$x\in S^*$ has an open neighborhood homeomorphic to an open ball in
$\Bbb R^{[\frac{n+1}2]}$, since all other type of singular points appear in
any such an open neighborhood, here a maximal singular point is either a
fixed point or a circle orbit.

Let $H$ denote the isotropy group of $x$. Assume that $n=2m$ or $2m+1$.
By the slice theorem, the orbit $x$ has an open neighborhood of the type,
$T^{[\frac{n-1}2]}\times _H D^{2m}/T^{[\frac{n-1}2]}=D^{2m}/H$, where
$D^{2m}$ is the open unit disk in $\Bbb R^{2m}$, and $H\cong T^{m-1}$.
Since $H$ acts linearly on $D^{2m}$ which also satisfies
(3.1.1) and (3.1.2), by Lemma 4.2 it follows that
$\partial D^{2m}(r)/H$ is a sphere for all $0<r<1$ and thus
$D^{2m}/H$ is homeomorphic to a disk.
\qed\enddemo

\vskip2mm

To prove that $M^*$ is homeomorphic to a sphere, we also need
the following two lemmas.

For a compact Lie group $G$, let $EG$ denote the universal principle
$G$-bundle. For a $G$-space $X$, the $i$-th equivariant cohomology
group $H_G^i(X):=H^i(X_G)$, where $X_G=EG\times_GX$ (cf. [Br]).

\vskip2mm

\proclaim{Lemma 4.4}

Let $X$ be a precompact $G$-space with $G=T^{[\frac{n-1}2]}$, and
let $S$ denote the singular set. If all isotropy groups are
connected, then

\noindent (4.4.1) For all $i\ge 0$, there is an isomorphism,
$H^i_G(X,S)\cong H^i(X^*,S^*)$.

\noindent (4.4.2) $H^i(X,S)=0$ (resp. $H_i(X, S)=0$) for all $i\le l$
if and only if $H^i_G(X,S)=0$ (resp. $H_{i, G}(X, S)=0$) for all $i\le l$.
\endproclaim

\vskip2mm

\demo{Proof} (4.4.1) Consider the natural map, $f: X_G \to X^*$. Note
that for $x\in (X^*- S^*)$, the preimage of $f^{-1}(x)$ is the
classifying space of the isotropy group at $x$ and thus contractible
since the isotropy group is trivial. Then (4..4.1) follows from
the proof of Proposition VII.1.1 of [Br] and the Vietoris-Begle
mapping Theorem.

(4.4.2) By a direct calculation.
\qed\enddemo

\vskip2mm

\demo{Proof of Theorem 4.1}

Since $M^*$ is a topological manifold of dimension at least
$5$ (Lemma 4.3), by Theorem 1.2
it suffices to show that $M^*$ is a homotopy sphere. Since
$M^*$ is simply connected (Lemma 1.5), $M^*$ is a homotopy sphere
if and only if $M^*$ is a homology sphere (Hurewicz theorem).
By the Poincar\'e duality, $M^*$ is a homology sphere if
$H_i(M^*)=0$ for all $i\le[\frac{m+1}2]$, where $n=2m$ or $(2m+1)$.

We now assume that $n\ge 8$. Note that by Theorem 3.6, there is a unique
$1$-simplex $\Delta_1^*$ of $S^*$ such that $\Delta_1^*\cap N^*=\emptyset$
(note that $\Delta^*_1$ is spanned by the two vertices of $S^*$ which are
not in $N^*$). By the transversality, $H_i(M^*)\cong H_i
(M^*-\Delta^*_1)\cong H_i((M^*-\Delta^*_1),N^*)$ for $1\le i\le m-2$
since $N^*$ is contractible. Observe that the inclusion
$N^*\subset (S^*-\Delta _1^*)$ (resp. $N\subset (S-p^{-1}(\Delta _1^*))$)
is a homotopy equivalence. Therefore, by Lemma 4.4,
$H_i((M^*-\Delta^*_1),
N^*)=0$ if and only if $H_i((M-p^{-1}(\Delta^*_1)),N)=0$ for
$i\le n-6$. 
Again by the transversality and Theorem 1.8,
$H_i((M-p^{-1}(\Delta^*_1)),N)\cong H_i(M,N)=0$ for $i\le n-6$.
Since $n\ge 8$, $n-6\ge [\frac{m+1}2]$ and thus the desired
result follows.
\qed\enddemo

\vskip 2mm

We conclude the section with the following application of Theorem 4.1.

\proclaim{Corollary 4.5}

Let $M$ be as in Theorem 4.1 of dimension $n\ge 8$. Assume that $n$ is
odd or $n$ is even and $\chi (M)=2$ (see Theorem 1.11). Then there is a
totally geodesic sphere $S^{n-4}\subset M$ and an embedded $3$-sphere
$S^3$ in $M$ such that $(M-S^3)$ is homotopy equivalent to $S^{n-4}$.
\endproclaim

\demo{Proof} Let $N$ and $\Delta^*_1$ be as in the proof of
Theorem 4.1. By Theorem 3.6 $N\simeq S^{n-4}$ is totally geodesic
and $N^*$ is contractible. Observe that $p^{-1}(\Delta^*_1)$ is homeomorphic
to $S^3$. Indeed, $p^{-1}(\Delta^*_1)$ is a totally 
geodesic $3$-sphere if $n$ is odd, and may not be totally geodesic
but contained in a totally geodesic $4$-sphere $p^{-1}(\Sigma (\Delta^*_1))$
if $n$ is even.  By Theorem 4.1,
$(M^*-\Delta^*_1)$ is contractible, and thus by Lemma 4.4
$$H_i((M-p^{-1}(\Delta^*_1)),(S-p^{-1}(\Delta^*_1)))\cong
H_i((M^*-\Delta^*_1),(S^*-\Delta^*_1))
\cong H_i((M^*-\Delta^*_1),N^*)=0$$
for all $i$, since $(S^*-\Delta^*_1)$ homotopy retracts to $N^*$. Therefore the inclusion
$j: (S-p^{-1}(\Delta^*_1))\longrightarrow (M-p^{-1}(\Delta^*_1))$
is a homology equivalence. By transversality the complement
$(M-p^{-1}(\Delta^*_1))$ is simply connected. Hence, by the Hurewicz
theorem, $j$ is a homotopy equivalence. The desired result follows.
\qed\enddemo

\vskip4mm

\head 5. Proof of Theorem B
\endhead

\vskip4mm

As explained in the introduction, in the proof of Theorem B
the main technical result is the following.

\vskip2mm

\proclaim{Theorem 5.1}

Let the assumptions be as in Theorem B. Then $M$ has a
submanifold $N$ such that

\noindent (5.1.1) $N$ is homeomorphic to $\Bbb CP^{n-1}$.

\noindent (5.1.2) The inclusion $N\hookrightarrow M$ is
$(2n-6)$-connected.
\endproclaim

\vskip2mm

\demo{Proof of Theorem B by assuming Theorem 5.1}

First, the Euler characteristic number of $M$ is $(n+1)$
(Corollary 1.12) and thus $M$ is homotopically
equivalent to $\Bbb CP^n$ (Theorem 1.10). By now Theorem
B follows from Theorem 5.1 and Corollary 1.4.
\qed\enddemo

\vskip2mm

\proclaim{Lemma 5.2}

Let $M$ be a closed $2n$-manifold. Assume that $M$ admits
an effective $T^n$-action with $(n+1)$ isolated fixed points and
$(M^*, S^*)\cong
(\Delta_n, \partial \Delta _n)$ as stratified spaces. Let $p: M\to M^*$ be
the orbit projection. If every face $\Delta_{n-1}\subset \partial
\Delta_n$, $p^{-1}(\Delta_{n-1})$ is diffeomorphic to
$\Bbb CP^{n-1}$, then $M$ is homeomorphic to $\Bbb CP^n$.
\endproclaim

\vskip2mm

Note that no curvature assumption is required in Lemma 5.2.

\vskip2mm

\demo{Proof} First note that $M$ is simply connected, since $M^*$ is
contractible and $T^n$ acts on $M$ with non-empty fixed point set.

We identify $M^*$ with $\Delta_n$. Choose a face
$\Delta _{n-1}\subset \Delta _n$ and identify $\Bbb CP^{n-1}$ with
$p^{-1}(\Delta_{n-1})$. Let $D(\eta)$ be a tubular neighborhood of
$\Bbb CP^{n-1}$, which is a closed $D^2$-bundle on $\Bbb CP^{n-1}$.
Let $W$ be the closure of the complement $M-D(\eta)$. Note that
$M=D(\eta ) \cup _{\partial W}  W$. Clearly, $W$ is homotopy equivalent
to $(M-\Bbb CP^{n-1})$.

We claim that $W$ is homologous equivalent to a point. Assuming this,
by duality the boundary $\partial W$ is a homology sphere. Since
$\partial W=S(\eta)$, the circle bundle of $D(\eta)$, by Gysin exact
sequence the Euler class of $\eta$ has to be a generator of $\Bbb CP^{n-1}$.
Therefore $\partial W=S^{2n-1}$. By the Van-Kampen theorem $\pi _1(W)=0$
since $\pi _1(M)=0$. Thus $W$ is diffeomorphic to the ball $D^{2n}$ and
so $M=D(\eta ) \cup _{S^{2n-1}} D^{2n}$ is unique up to homeomorphism.

It remains to prove $(M-\Bbb CP^{n-1})$ is homologous to a point.
By the assumption $S\cong S(\Bbb CP^{n},T^n)$. In other words, it 
is union of $n$ copies of $\Bbb CP^{n-1}$ such that every two 
copies intersects in a $\Bbb CP^{n-2}$, etc. Observe that
$(S-\Bbb CP^{n-1})$ is contractible (it obviously
deformation retracts to the unique fixed point outside $\Bbb CP^{n-1}$).
Applying Lemma 4.4 to $(M-\Bbb CP^{n-1})$ we get
$$H_i((M-\Bbb CP^{n-1}),
(S-\Bbb CP^{n-1}))\cong H_i((\Delta _n- \Delta_{n-1}), (\partial
\Delta _n- \Delta_{n-1}))=0$$
for all $i$. Therefore $H_i(M-\Bbb CP^{n-1})=0$ for all $i\ge 1$. The
desired result follows.
\qed\enddemo

\vskip2mm

\demo{Proof of Theorem 5.1}

(5.1.1) By the discussion in Section 2, we may assume
(3.1.1) and (3.1.2). The singular set $S^*$ is isomorphic to
$\Delta_n^{n-2}$ (Theorem 3.6) and $M^*$ is homeomorphic to
$S^{n+1}$ (Theorem 4.1). These two facts are crucial for the
following construction of the desired $N$.

Let $v_1,...,v_n, v_{n+1}$ denote the vertices of $S^*$, and let
$A$ denote the simplicial subcomplex of $S^*$ which consists of
$(n-2)$-simplices which do not contain $v_{n+1}$. Then the total
space $|A|$ of $A$ is homeomorphic
to $S^{n-2}$ (Theorem 3.6). A priori, $|A|$ is a codimension $3$ knot
in $M^*$. By [St], $|A|$ must be a topological trivial knot and thus
bounds a $(n-1)$-ball $D\subset M^*$. We may assume that $D\cap S^*=A$
(since $(M^*-(S^*-|A|))$ is homeomorphic to $\Bbb R^{n+1}$).  Put $N=p^{-1}
(D)$. Note that $N$ is a closed invariant $(2n-2)$-manifold. Since
every top simplex of $S$ in $N$ is a totally geodesic $\Bbb CP^{n-2}$,
by Lemma 5.2 $N$ is homeomorphic to $\Bbb CP^{n-1}$.

(5.1.2) Note that $\Bbb CP^{n-2}\to  N\approx \Bbb CP^{n-1}$ is
$(2n-4)$-connected. Hence, by Theorem 1.8 the inclusion
$i: \Bbb CP^{n-2}\to  N\hookrightarrow M$ is
$(2n-6)$-connected.
\qed\enddemo

\vskip4mm

\head 6. Proof of Theorem A.
\endhead

\vskip4mm

First, by the discussion in Section 2, we may reduce the
proof of Theorem A to the case in (3.1.1) and (3.1.2).

The proof of Theorem A for $n=8$ and $9$ follows from Lemmas 6.1
and 6.2.

\vskip2mm

\proclaim{Lemma 6.1}

Let $M$ be as in Theorem A. If $n=8$, then $M$ is homeomorphic
to $S^8$, $\Bbb CP^4$ or $\Bbb HP^2$.
\endproclaim

\vskip2mm

\demo{Proof} We will show that $M$ is
homeomorphic to $S^8$, $\Bbb HP^2$ or $\Bbb CP^4$ corresponding
respectively to $\chi(M)=2, 3$ or $5$ (Corollary 1.12).

Case 1. Assume that $\chi(M)=2$.

We only need to show that $M$ is a homotopy sphere (Theorem 1.2).
Since $\chi(M)=2$, by Poincar\'e duality and the Hurewicz
theorem it is easy to see that $M$ is a homotopy sphere if
$M$ is $3$-connected, which follows trivially from Corollary 4.5 and
the Alexander duality.

Case 2. Assume that $\chi(M)=3$.

Let $S^4$ denote a totally geodesic sphere which is a $4$-simplex
of $S$. By the argument similar to the above, one may
verify that $(M-S^4)$ is contractible and so homeomorphic to the open
ball $\text{int}(D^8)$ (by the h-cobordism theorem),
and thus $M=D(\nu )\cup _\partial \bar D^8$, where $D(\nu )$
is the normal disk bundle over $S^4$ with fiber $D^4$.
Since $S^4$ is the fixed point component of a circle isotropy
group of $T^3$-action without any finite order isotropy group,
the normal bundle
$\nu$ must be a complex $2$-bundle with normal sphere
bundle $S^7$. This shows that the second Chern class
of the bundle is a generator of $H^4(S^4,\Bbb Z)$ and
$\nu$ is equivalent to the normal bundle of $\Bbb HP^1
\subset \Bbb HP^2$. In particular,
this gives a diffeomorphism $f: (M-\text{int}(D^ 8)) \to (\Bbb HP^2 -
\text{int}(D^ 8))$. Clearly the restriction of $f$ on the boundary
may be extended to a homeomorphism of $D^8$ by radical extension
(indeed one may prove the diffeomorphism type is unique using
some deep topology in this case). The desired result follows.

Case 3. Assume $\chi(M)=5$.

By Theorem B, it suffices to prove that $M$ is homotopy equivalent to
$\Bbb CP^4$. To achieve this, we only need to prove that the cohomology
ring of $M$ is the same as that of $\Bbb CP^4$.
Let $v_1, v_2, v_3, v_4, v_5$ be the vertices of $S^*\cong
\Delta _4^2$ (Theorem 3.6). Recall that a face $S^*$, saying $[v_3v_4v_5]$,
corresponds to a totally geodesic $\Bbb CP^2$ in $M$. By removing a face
spanned by $v_3, v_4, v_5$ in $\Delta _4^2$, we claim that the complement
$M-\Bbb CP^2$ is homotopy equivalent to $S^2$, which corresponds to a
segment $[v_1v_2]$ with vertices $v_1, v_2$. We omit the details since
the argument is similar to the proof of Corollary 4.5.  Therefore, by
Alexander duality Theorem we get that $H_i(M)\cong H_i(\Bbb CP^2)$ for
$i\le 5$. By Poincare duality it suffices to prove that the self
intersection
number of the cycle $[\Bbb CP^2]$ is $\pm 1$. This is because that the
simplex of $S$ represented by the face $[v_1v_2v_3]$ intersects with
$\Bbb CP^2$ at a single point $v_3$. The desired result follows.
\qed\enddemo

\vskip2mm

\proclaim{Lemma 6.2}

Let $M$ be as in Theorem A. If $n=9$, then $M$ is
homeomorphic to a sphere.
\endproclaim

\vskip2mm

\demo{Proof} Since $S$ has a $5$-simplex (Theorem 3.6),
$M$ is $3$-connected (Theorem 1.8). By Poincar\'e duality,
it suffices to prove that $H_4(M)=0$. This follows easily
from Corollary 4.5 and the Alexander duality.
\qed\enddemo

\vskip4mm

\head 7. Proof of Theorem D For $n=6$
\endhead

\vskip4mm

The goal of this section is to establish the following result.

\vskip2mm

\proclaim{Proposition 7.1}

Let $M$ be a closed simply connected $6$-manifold. Assume
$M$ admits a $T^2$-action satisfying following conditions:

\noindent (7.1.1) The $T^2$-fixed point set is non-empty and
there is no circle subgroup with fixed point set codimension $2$.

\noindent (7.1.2) All isotropy groups are connected.

\noindent (7.1.3) Each fixed point component of any isotropy group
is either a $2$-sphere or a point.

\noindent Then $M$ is diffeomorphic to $S^6$ if and only if
$b_2(M)=0$.
\endproclaim

\vskip2mm

\remark{Remark \rm 7.2} Proposition 7.1 will be true without
the restriction (7.1.2). A reason we do not
include a proof for this more general situation is
to avoid some unnecessary technical complexity (compare to [Ro1]);
in our circumstance we do have (7.1.2) (see (3.1.1)).
\endremark

\vskip2mm

An interesting consequence of Proposition 7.1 is that any $T^2$-action
on $S^3\times S^3$ does not satisfy (7.1.2) or (7.1.3).

\vskip2mm

\demo{Proof of Theorem D for $n=6$} By Lemma 2.1, we may assume
(7.1.1) and (7.1.2). Since a fixed point component is an
orientable  totally geodesic submanifold of even-codimension,
(7.1.3) follows from the positive curvature assumption. By
now Theorem D for $n=6$ follows from Proposition 7.1.
\qed\enddemo

\vskip2mm

\proclaim{Lemma 7.3}

Let $M$ be as in Proposition 7.1. If $b_2(M)=0$, then the
Euler characteristic $\chi(M)=2$.
\endproclaim

\vskip2mm

\demo{Proof} Since $b_2(M)=0$, $\chi(M)=2-b_3(M)\le 2$.
By (7.1.1), from the isotropy representation of the
$T^2$-action at an isolated fixed point one concludes
that there is a circle subgroup whose fixed point
set contains a sphere $S^2$ (see (7.1.3)). This implies that
$\chi (M)\ge 2$ (see Theorem 1.6) and therefore
$\chi (M)=2$ and $b_3(M)=0$.
\qed\enddemo

\vskip2mm

\proclaim{Lemma 7.4}

Let $M$ be as in Proposition 7.1. Then $H_2(M)$ has no torsion.
\endproclaim

\vskip2mm

\demo{Proof} Note that $H_2(M)\cong \pi_2(M)$ (Hurewicz theorem).
Since the singular set $S$ is a finite union of at most
$2$-dimensional submanifolds, by transversality $\pi_2(M)
\cong \pi_2(M-S)$. Thus it suffices to show that
$\pi_2(M-S)$ has no torsion.

By the homotopy exact sequence of the
fibration, $T^2\to (M-S)\to (M^*-S^*)$,
$$0\to\pi_2(M-S)\to \pi_2(M^*-S^*)\to \Bbb Z^2\to 0$$
it reduces to show that $\pi_2(M^*-S^*)$ has no torsion.
By (7.1.1) and (7.1.2), it is easy to see that $M^*$ is a
topological $4$-manifold. Since $M^*$ is simply
connected (see Lemma 1.5), by transversality $(M^*-S^*)$ is
simply connected and thus
$\pi_2(M^*-S^*)\cong H_2(M^*-S^*)$
(Hurewicz theorem). By Alexander duality,
$H_2(M^*-S^*)\cong H^2(M^*, S^*)$. Since both $H^1(S^*)$ and
$H^2(M^*)$ are torsion free, the desired result follows by the
long exact sequence of the pair $(M^*, S^*)$.
\qed\enddemo

\vskip2mm

\demo{Proof of Proposition 7.1}

If $b_2(M)=0$, by Lemma 7.4 we know that $\pi_2(M)=0$ and by
Lemma 7.3 $b_3(M)=0$. By duality and Hurewicz theorem this implies
that $M$ is a homotopy sphere. By Theorem 1.2 $M$ is homeomorphic
and thus diffeomorphic to $S^6$.
\qed\enddemo

\vskip4mm

\head 8. Proof of Theorem D For $n=7$
\endhead

\vskip4mm

The goal of this section is to establish the following topological
result.

\vskip2mm

\proclaim{Proposition 8.1}

Let $M$ be a closed simply connected $7$-manifold.
Assume that $M$ admits a $T^3$-action such that

\noindent (8.1.1) There is no circle subgroup with fixed point set of
codimension $2$.

\noindent (8.1.2) All isotropy groups are connected.

\noindent (8.1.3) Each fixed point component of any isotropy group
is either a circle or a lens space of dimension $3$.

\noindent Then $M$ is homeomorphic to $S^7$ if and only if
$b_2(M)=0$.
\endproclaim

\vskip2mm

Note that Proposition 8.1 will still be true without (8.1.2);
compare to Remark 7.2.

\vskip2mm

\demo{Proof of Theorem D for $n=7$} By Lemma 2.1, we may assume
(8.1.1) and (8.1.2). Since any fixed point component of an
isotropy group is a closed totally geodesic manifold with
the maximal symmetry rank, (8.1.3) holds (see Theorem 1.9).
By now Theorem D for $n=7$ follows from Proposition 8.1.
\qed\enddemo

\vskip2mm

We first determine the singular structure and the 
homeomorphic type of the orbit space.

\vskip2mm

\proclaim{Lemma 8.2}

Let the assumptions be as in
Proposition 8.1. Then $S^*$ is isomorphic to $\Delta_3^1$
and $M^*$ is homeomorphic to $S^4$.
\endproclaim

\vskip2mm

\demo{Proof} By Theorem 4.1, $M^*$ is a closed topological
$4$-manifold. By Lemma 1.5 $M^*$ is simply connected. We further claim
that $M^*$ is a homotopy $4$-sphere. In fact, by [FR1] Lemma 5.2
we see that $H_2(M/H;\Bbb Q)=0$, where $H\cong T^2$ is an isotropy group.
Observe that $T^3/H$ acts on $M/H$ with non-empty fixed point set. Since
$M^*=(M/H)/(T^3/H)$, by [FR1] Lemma 5.2 once again we get that
$H_2(M^*;\Bbb Q)=0$. By Poincare duality this implies that $M^*$ is a
homotopy $4$-sphere.

Consider the homotopy exact sequence of $T^3\to (M-S)\to (M^*-S^*)$,
$$0\to \pi _2(M-S)\to \pi _2(M^*-S^*)\to \Bbb Z^3\to 0.$$
By the transversality (see (8.1.1)),
$\pi _2(M-S)\cong\pi _2(M)=0$ and thus
$\pi _2(M^*-S^*)\cong\Bbb Z^3$.
By Alexander duality $H_1(S^*)
\cong H_2(M^*, S^*)\cong H^2(M^*-S^*)\cong \Bbb Z^3$.

Note that the local singular structure described in Lemma 3.8
remains valid in our situation: around each isolated circle
orbit, there are exactly three $3$-dimensional strata sharing
the circle orbit. Hence, at each vertex of $S^*$ there are
exactly three edges going out from the vertex. Since $S^*$
is a graph satisfying $b_1(S^*)=3$, it is easy to see that
$S^*\cong \Delta _3^1$.
\qed\enddemo

\vskip2mm

\demo{Proof of Proposition 8.1}

We claim that any $3$-simplex of $S$ is homeomorphic to $S^3$
and $(M-S^3)$ is homology equivalent to a three sphere.
By Alexander duality, it follows easily that $M$ is a homology 
sphere  and thus a homotopy sphere since $M$ is simply connected.
By Theorem 1.2 $M$ is homeomorphic to a sphere.

Let $L$ denote a $3$-simplex of $S$. By (8.1.3), $L$ is a lens space.
By Lemma 8.2, there is a unique $3$-simplex $L'$ of $S$ such that
$L^*\cap L^{'*}=\emptyset$. By (8.1.2) and Lemma 4.4,
$H_i((M-L'),(S-L'))\cong H_i((M^*-L^{'*}),(S^*-L^{'*}))=0$
for all $i$ (see Lemma 8.2). Consequently, $(M-L')$ is homology equivalent
to $(S-L')$ which, in turn, is homology equivalent to $L$ (see Lemma 8.2).
On the other hand, by (8.1.1) and the transverality $\pi_1(M-L')\cong
\pi_1(M)=1$. Hence, $H_1(L)=0$ and so $L$ is simply connected and therefore
homeomorphic to $S^3$. Likewise, $L'$ is also homeomorphic
to $S^3$. The desired result follows.
\qed\enddemo


\vskip4mm

\Refs
\nofrills{References}
\widestnumber\key{APS11}

\vskip3mm

\ref
\key AW
\by S. Aloff; N. R. Wallach
\pages 93-97
\paper An infinite family of $7$-manifolds admitting positive curved
Riemannian structures
\jour Bull. Amer. Math. Soc.
\vol 81
\yr 1975
\endref

\ref
\key Br
\by G. Bredon
\paper Introduction to compact transformation groups
\jour Academic Press
\vol 48
\yr 1972
\endref


\ref
\key Es
\by J.-H Eschenburg
\pages 469-480
\paper New examples of manifolds with strictly positive curvature
\jour Invent. Math
\yr 1982
\vol 66
\endref

\ref
\key FMR
\by F. Fang; S. Mendonca; X. Rong
\pages
\paper A connectedness principle in the geometry of positive
curvature
\jour Preprint
\yr 2001
\vol
\endref

\ref
\key FR1
\by F. Fang; X. Rong
\pages 641-674
\paper Positive pinching, volume and second Betti number
\jour Geom. Funct. Anal.
\yr 1999
\vol 9
\endref

\ref
\key FR2
\by F. Fang; X. Rong
\pages
\paper Positively curved manifolds of maximal discrete
symmetry rank
\jour Amer. J. Math 
\yr To appear
\vol
\endref

\ref
\key FR3
\by F. Fang; X. Rong
\pages
\paper Positively curved manifolds of almost maximal
symmetry rank in dimensions $6$ and $7$
\jour In preparation
\yr
\vol
\endref

\ref
\key Fr
\by M. Freedman
\pages 357-453
\paper Topology of Four Manifolds
\jour J. of Diff. Geom.
\yr 1982
\vol 28
\endref

\ref
\key Gro
\by K. Grove
\pages 31-53
\paper Geometry of, and via symmetries
\jour Univ. Lecture Ser., Amer. Math. Soc., Providence, RT
\yr 2002
\vol 27
\endref

\ref
\key GKS
\by S. Goette; N. Kitchlook; K. Shankar
\pages
\paper Diffeomorphism type of the Berger space $SO(5)/SO(3)$.
\jour Preprint
\yr
\vol
\endref


\ref
\key GS
\by K. Grove, C. Searle
\pages 137-142
\paper Positively curved manifolds with maximal symmetry-rank
\jour J. Pure Appl. Alg
\yr 1994
\vol 91
\endref

\ref
\key GZ
\by K. Grove; W. Ziller
\pages 331-367
\paper Curvature and symmetry of Milnor spheres
\jour Ann. of Math
\vol 152
\yr 2000
\endref

\ref
\key Hs
\by W. Hsiang
\pages
\paper Cohomology theory of topological transformation groups
\jour Ergebnisse der Mathematik und inere Grenzgebiete
\yr 1975
\vol 85
\endref

\ref
\key HK
\by W. Hsiang, B. Kleiner
\pages 615-621
\paper On the topology of positively curved $4$-manifolds
with symmetry
\jour J. Diff. Geom
\yr 1989
\vol 30
\endref

\ref
\key Ko
\by S. Kobayashi
\pages
\paper Transformation groups in differential geometry
\jour Springer-Verlag Berlin Heidelberg New York
\vol
\yr 1972
\endref

\ref
\key PS
\by T. P\"uttmann; C. Searle
\pages 163-166
\paper The Hopf conjecture for manifolds with low cohomogeneity
or high symmetry rank
\jour Proc. Amer. Math. Soc.
\vol 130
\yr 2002
\endref

\ref
\key Ro1
\by X. Rong
\pages 157-182
\paper Positively curved manifolds with almost maximal symmetry rank
\jour Geometriae Dedicata
\yr 2002
\vol 59
\endref

\ref
\key Ro2
\by  X. Rong
\pages 397-411
\paper On the fundamental groups of compact manifolds of positive
sectional curvature
\jour Ann. of Math
\yr 1996
\vol 143
\endref

\ref
\key Sm
\by S. Smale
\pages 391-466
\paper Generalized Poincar\'e conjecture in dimension $>4$
\jour Ann. of Math
\yr 1961
\vol 74
\endref

\ref
\key St
\by J. Stallings
\pages  490-503
\paper On topologically unknotted spheres
\jour Ann. of Math
\yr 1963
\vol 77
\endref

\ref
\key Su
\by D. Sullivan
\pages
\paper Triangulating homotopy equivalences and homeomorphisms
\jour   Geometric Topology Seminar Notes, in ``The Hauptvermutung
Book'' , A collection of papers on the topology of manifolds,
K-monographs in Math.
\yr edited by A. Ranicki, Kluwer Academic Publishers, 1995
\vol 1
\endref

\ref
\key Wi
\by  B. Wilking
\pages
\paper Torus actions on manifolds of positive sectional curvature
\jour preprint (August 2002)
\yr
\vol
\endref

\ref
\key Ya
\by D. Yang
\pages 531-545
\paper On the topology of nonnegatively curved simply connected
$4$-manifolds with discrete symmetry
\jour Duke Math. J.
\yr 1994
\vol 74
\endref

\endRefs
\enddocument